\documentclass[a4, 10pt,leqno]{amsart}

\usepackage{mathabx}

\usepackage{a4wide}

\usepackage[foot]{amsaddr}
\usepackage[active]{srcltx}
\usepackage[utf8]{inputenc}
\usepackage[T1]{fontenc}
\usepackage{upgreek}
\DeclareMathAlphabet{\mathdutchcal}{U}{dutchcal}{m}{n}
\SetMathAlphabet{\mathdutchcal}{bold}{U}{dutchcal}{b}{n}
\DeclareMathAlphabet{\mathdutchbcal}{U}{dutchcal}{b}{n}
\DeclareSymbolFont{myletters}{OML}{ztmcm}{m}{it}
\DeclareMathSymbol{\nicelambda}{\mathord}{myletters}{"15}
\usepackage{amsmath,amssymb,amscd,amsthm,graphicx,enumerate}
\usepackage{mathrsfs}
\usepackage{amssymb,amsthm,a4wide}
\usepackage{nicefrac}
\newcounter{example}
\newenvironment{example}[1]{\refstepcounter{example}\par\medskip
	\noindent \textsc{\small Example~\theexample. #1} \rmfamily\hspace{-2pt}}{\medskip}

\def\bysame{\leavevmode\hbox to3em{\hrulefill}\thinspace}
	
\usepackage{pifont}

\usepackage{color}
\usepackage[usenames,dvipsnames]{xcolor}
\usepackage{xcolor}
\usepackage{hyperref}
\hypersetup{
	linktoc=page,
	colorlinks,
	linkcolor={Maroon},
	citecolor={Maroon},
	urlcolor={Maroon}
}
\hypersetup{breaklinks=true}
\allowdisplaybreaks

\usepackage{enumerate}
\makeatletter
\@namedef{subjclassname@2020}{%
	\textup{2020} Mathematics Subject Classification}
\makeatother

\usepackage[a4paper,bindingoffset=0in,%
left=1.02in,right=1.02in,top=1.22in,bottom=1.22 in,%
footskip=1mm]{geometry}

\theoremstyle{plain}

\newtheorem*{theorem*}{Theorem}

\newtheorem{lemma}{Lemma}[section]
\newtheorem{theorem}[lemma]{Theorem}
\newtheorem{proposition}[lemma]{Proposition}
\newtheorem{corollary}[lemma]{Corollary}

\theoremstyle{definition}
\newtheorem{definition}[lemma]{Definition}
\newtheorem{remark}[lemma]{Remark}

\theoremstyle{remark}

\numberwithin{equation}{section}

\newcommand{\R}{\mathbb R}
\newcommand{\Z}{\mathbb Z}

\newcommand{\N}{\mathcal N}
\newcommand{\K}{\mathcal K}

\newcommand{\G}{\mathcal G}

\renewcommand{\setminus}{\smallsetminus}

\DeclareMathOperator{\supp}{supp}

\makeatletter
\newcommand{\tpmod}[1]{{\@displayfalse\pmod{#1}}}
\makeatother

\DeclareMathOperator{\Ima}{Im}
\newcommand{\sph}{\mathbb{S}}

\newcommand{\transv}{\mathrel{\text{\tpitchfork}}}
\makeatletter
\newcommand{\tpitchfork}{%
	\vbox{
		\baselineskip\z@skip
		\lineskip-.52ex
		\lineskiplimit\maxdimen
		\m@th
		\ialign{##\crcr\hidewidth\smash{$-$}\hidewidth\crcr$\pitchfork$\crcr}
	}%
}
\makeatother

\newcommand{\M}{\mathcal{M}}
\newcommand{\TM}{\mathcal{TM}}
\newcommand{\E}{\mathcal{E}}

\begin{document}

\title[\scriptsize {Rigidity of ALE vector bundles}]{\small On rigidity of ALE vector bundles} 
\author[{\protect \scriptsize F. Asadi}]{Fatemeh Asadi$^1$}
\address{$^1$Department of Mathematics, Azarbaijan Shahid Madani University, P.O. Box 53714-161, Tabriz, Iran.}
\email{\href{mailto:f_asadi@azaruniv.ac.ir}{f\_asadi@azaruniv.ac.ir}}

\author[\protect\scriptsize Z. Fathi]{Zohreh Fathi$^2$}
\email{\href{mailto:z.fathi@ipm.ir}{z.fathi@ipm.ir}}

\address{$^{2,3}$School of Mathematics\\
	Institute for Research in Fundamental Sciences (IPM) \\
	P.O. Box 19395-5746\\
	Tehran\\
	Iran
}
\email{\href{mailto:slakzian@ipm.ir}{slakzian@ipm.ir}}

\author[\protect \scriptsize S. Lakzian]{Sajjad Lakzian$^{3^*}$}
\address{$^{3}$Department of Mathematical Sciences\\ Isfahan University of  Technology (IUT) \\ Isfahan 8415683111, Iran}

\email{\href{mailto:slakzian@iut.ac.ir}{slakzian@iut.ac.ir}}

%

\subjclass[2020]{Primary: 53Cxx, 57Rxx; Secondary: 57Nxx, 58Axx}
\keywords{asymptotically conical manifold, ALE manifold, bundles, asymptotically flat manifolds, flat manifolds, spherical space forms}
\thanks{*\textit{the corresponding author}}
\maketitle
%
%
%
%

\begin{abstract}
	\textsl{We discuss topological rigidity of vector bundles with asymptotically conical ({\sf AC}) total spaces of rank $\ge 2$ with a sufficiently connected link; our focus will mainly be on {\sf ALE} (asymptotically locally Euclidean) bundles. Within the smooth category, we topologically classify all {\sf ALE} tangent bundles by showing only $\sph^2$, $\mathbb{RP}^2$ and open contractible manifolds admit {\sf ALE} tangent bundles. We also discuss other interesting topological and geometric rigidities of {\sf ALE} vector bundles.}
\end{abstract}
\date{\today}
%
%
%
%
%

\section{Introduction}
\par Asymptotically conical ({\sf AC}) manifolds are in general expected to model ``reasonable'' open manifolds $\M$ with asymptotically no slower than quadratic curvature decay that also admit unique tangent metric cones at infinity. For example, this is true under asymptotically faster than quadratic curvature decay (a.k.a. weak asymptotic flatness) and the unique tangent metric cone property~\cite{PT}. The latter conditions alone do not impose any topological restrictions on $\M$ however with the stronger variants of the said conditions, with nonnegative Ricci curvature assumptions and with enough volume growth, topological rigidity on $\M$ emerges toward becoming contractible. For instance, strong asymptotic flatness plus strong holonomy control implies ends of $\M$ are {\sf AC} with links that are fibrations over spherical space forms~\cite{CL}, while faster than quadratic curvature decay along with Euclidean volume growth implies $\M$ is asymptotically locally Euclidean ({\sf ALE})~\cite{BKN}; in the latter case, if in addition $\M$ is Ricci nonnegatively curved, one can bound the size of its fundamental group by the reciprocal of asymptotic volume growth~\cite{Li,Anders} and if the growth is past Prelman's threshold, the manifold must be contractible~\cite{Per}.
\par {\sf ALE} manifolds form a pivotal special case of {\sf AC} manifolds that organically arise in several fields; they arise as gravitational instantons~\cite{Haw,EH1,GH} as well as singularity models for 4D Ricci flow~\cite{BZ,Sim}. A special case, is {\sf AE} manifolds that are main objects of study as spatial models in general relativity.
\par With extra structures, partial classifications of {\sf AC} manifolds are possible. For example, in the K\"{a}hler setting, there is a more rigidity to {\sf AC} manifolds; indeed, it is shown in~\cite{HRS}, that all K\"{a}hler {\sf ALE} manifolds are obtained as resolutions of deformations of an orbifold singularity~$\nicefrac{\mathbb{C}}{G}$. Also all Clalabi-Yau {\sf AC} manifolds polynomially asymptotic (up to diffeomorphism) to the hyperquadric in $\mathbb{C}^{n+1}$ are classified~\cite{CH}. 
\par The extra structure we are interested in, is that of a vector bundle. Total spaces of vector bundles are better behaved topologically and posses more structure than their base spaces, a feature that makes them and specially tangent bundles (recall a tangent bundle can be characterized by existence of a globally integrable tangential structure) suited to host many important constructions especially in the context of {\sf AC} manifolds. To name a few important instances, the Eguchi-Hanson metric, and one of the Candelas and de la Ossa's metrics in 3D along with the Stenzel's metrics~\cite{EH1,Cand-Dela,Stenzel} (the latter generalizes the other two) are all constructions on the (co)tangent bundle of spheres while another important metric due to Candelas and de la Ossa is on $\mathcal{O}_{\mathbb{P}^1}(-1)^{\oplus 2}$~\cite{Cand-Dela}. Indeed, in complex dimensions higher than 3, Stenzel's construction provide the only Ricci flat Calabi-Yau ALE manifolds that are polynomially asymptotic to the hyperquadric as is recently shown in \cite{CH}.
\par Our goal in these notes is to explore the mostly topological rigidity of {\sf AC} vector bundles with a focus on {\sf ALE} bundles, especially tangent bundles as prototypes of vector bundles and of open manifolds admitting an almost complex structure; and as an application, to use these topological rigidities to also derive some geometric rigidity properties concerning the Riemannian metrics. One main aspect of these notes is the topological classification of {\sf ALE} tangent bundles and some non-existence theorems on {\sf ALE} bundles. 
\subsection{Set up and Notations}
\par Throughout these notes, $\widetilde{Z}$ denotes the universal cover of $Z$.  All spaces arising are smooth manifolds unless otherwise stated and the groups acting on them are discrete freely acting subgroups of diffeomorphisms. 
\subsubsection{Asymptotically Conical Manifolds}We will consider the following asymptotically conical spaces. 
\begin{definition}[Asymptotically Conical Manifolds]
	\hfill\\
	\begin{enumerate}
		\item  \textbf{\small \textsf{top-AC:}} A manifold $\M^n$ is said to be topologically asymptotically conical ({\sf top-AC}), if there exists a compact subset $\K \Subset \M$ and manifolds $\N^{n-1}_j$, $j=1, \cdots, m$, such that $\M \smallsetminus \K$ is the disjoint union of ends $\M_j$, each \emph{homeomorphic} -- via the restriction of a homeomorphism $\varphi$ -- to a punctured cone with link $\N_j$.  In particular, for any $p_j\in \M_j$, one assigns the coordinate $x_j$ defined via $\varphi^{-1}\left(p_j\right)=(x_j,q_j)$ where $q_j \in \N_j$. We further assume the increase in $x_j$ indicates getting further along the end; meaning if $\K \Subset U \subset  \K'$ for some open set $U$ and compact set $\K'$, then the values of $x_j$ on $\M\setminus \K'$ are strictly larger than those on $\M \setminus \K$. 
		\item  \textbf{\small \textsf{diff-AC:}} $\M$ from item (1) is said to be {\sf diff-AC} if ends are \emph{diffeomorphic}, via the restriction of a diffeomorphism $\varphi$, to the aforementioned cones.
		\item  \textbf{\small \textsf{Riem-AC:}} Suppose $\M$ is a {\sf diff-AC} and the metrics $g_{_{\N_j}}$ are given on $\N_j$ with the corresponding cone metrics $\widehat{g}_j := dx_j^2 + x_j^2 g_{_{\N_j}} $ (so, $x_j$ is now the distance to the vertex). The Riemannian manifold $\left(\M, g\right)$ is called a {\sf Riem-AC} manifold  if the restriction of the metric to the end $\M_j$, $g_j$, is of the form
		$
		\varphi^* g_j = \widehat{g}_j + o(1)
		$.
		Sometimes one may also require a fall off rate $\tau_j>0$ on $g_j$ and on derivatives up to order $d_j$ which means requiring the weighted Sobolev norm decay
		\[
		\sum_{\ell=0}^{d_j} x_j^{l} \|\nabla^\ell_{\widehat{g}_j}\left( \varphi^* g_j - \widehat{g}_j\right) \|_{\widebar{g}_j} = O(x_j^{-\tau_j}), \quad \text{as} \quad x_j \to \infty,
		\]
	\end{enumerate}
\end{definition}
\begin{definition}
	\hfill
	\begin{enumerate}
		\item
		A {\sf top-AC} manifold $ \M$ is called \emph{simple} whenever for all $j$, $\N_j= \nicefrac{\widetilde{\N}_j}{\Upgamma_j}$ for a finite freely acting subgroup $\Upgamma_j \le \mathrm{Diff}(\widetilde{\N_j})$. One can similarly define simple {\sf diff-AC} and {\sf Riem-AC} manifolds. For the latter, one assumes $\Upgamma_j \le \mathrm{Isom}(\widetilde{\N_j})$.
		\item A simple {\sf top-AC} ({\sf diff-AC} or {\sf Riem-AC}) manifold in which the groups $\Upgamma_j$'s are all orientation preserving, is called a \emph{proper} {\sf top-AC} ({\sf diff-AC} or {\sf Riem-AC}) manifold. 
	\end{enumerate}
\end{definition}
\par A particular case is when $\M$ is simple and $\widetilde{\N_j}$'s are $\sph^{n-1}$ (up to homeomorphism, diffeomorphism (nonexotic structure) or also equipped with standard metric); in this case, we get the definition of {\sf top-ALE} ({\sf diff-ALE} and {\sf Riem-ALE}) manifolds. Furthermore, allowing exotic spheres, one also defines exotic {\sf diff-ALE} manifolds. One can analogously define complex {\sf ALE} manifolds whose ends are diffeomorphic to quotients $\mathbb{C}^*$ by finite freely acting subgroups of $\mathrm{U}(n)$ (or $\mathrm{SU}(n)$ for proper ones). 
\par \emph{It is worth noting that, in many places in literature, an {\sf ALE} manifold refers to a \emph{proper} {\sf Riem-ALE} manifold.} 
\begin{definition}\label{defn:regular}
	A single ended {\sf top-AC} manifold $\M$ is called ``regular'' whenever the homeomorphism $\varphi$ extends to $\partial \K$; in particular, a {\sf Riem-AC} manifold can always assumed to be regular (perhaps by taking a larger $\K$).
\end{definition}
\begin{definition}
	When $\K$ is allowed to be a $\sigma$-compact and $\M\setminus \K$ to be comprised of countably of many connected components, $\M$ is called a $\sigma$-{\sf AC} manifold (with {\sf top}, {\sf diff} and {\sf Riem} variations).
\end{definition}
\par When $\M$ is {\sf top-AC} with one end, the fundamental group of the end $\M\setminus \K$ is $\Upgamma$ which coincides with $\pi_1(\N)$, this is also the fundamental group at infinity thus written as $\pi_1^\infty(\M)$.
\subsubsection{Bundles}\label{subsubsec:vb}
\par Throughout these notes, $\E^{r\ge 2}$ denotes vector bundles of rank $r\ge 2$ and we note the simple yet key fact that ``\emph{any vector bundle of rank $\ge 2$ has only one end, in particular $\TM^n$ has only one end when $n \ge 2$.}'' 
\par Also we assume $\M^n$ is a boundary-less manifold with $n\ge 2$, consequently so are the total spaces of all vector bundles $\mathcal{E}$ over $\M$.
\par Given a fiber bundle $\mathcal{F}$ -- with slight abuse of terminology -- we will write $\mathcal{F}$ for both the bundle and the total space of the bundle as a topological space. 
\par We work under the running assumption that $\E$ has a finite fundamental group. $\E_{\sf slit}$ will denote the slit bundle i.e. bundle minus the zero section. 
\subsubsection{Asymptotic Flatness and Holonomy Control Properties}
\begin{definition}
	Let $\left(\M, g \right)$ be a complete open manifold.
	\begin{enumerate}
		\item $\M$ is called ``weakly asymptotically flat'' if for some (thus all) $p\in \M$,
		$$
		\limsup_{r\to \infty} \left\{ \|\mathrm{Rm}\|(x), r^2  \right\} = 0,
		$$
		 where $r:= d(x,p)$ i.e. $\M$ asymptotically has faster than quadratic curvature decay;
		\item $\M$ is called ``strongly asymptotically flat'' if for any $p\in \M$, 
		$
		\|\mathrm{Rm}\|(x) \le \nicefrac{\Uplambda(r)}{r^2},
		$
		for a non-increasing positive function $\Uplambda$ with $\int_1^\infty \; \nicefrac{\Uplambda(s)}{s}\; ds < \infty$  i.e. $\M$ asymptotically has faster than quadratic curvature decay in a controlled way. 
	\end{enumerate}
	A main reference where these conditions have been studied in, is~\cite{Abesch}. 
\end{definition}
\begin{definition}[Strong holonomy control~\cite{CL}]
	Let $\Delta\Theta(\upgamma)$ denote the rotation along a geodesic loop $\upgamma$ (the angle between vectors and their parallel transport along $\upgamma$) in a Riemannian manifold $(\M,g)$. $(\M,g)$ satisfies strong holonomy control at scale $\kappa$ whenever there exists a compact set $\K$ and and a positive function $\varepsilon(r)$ with $\varepsilon(r)\downarrow 0$ as $r\downarrow 0$ such that 
	$
	\|	\Delta \Uptheta\| \le \varepsilon(l)
	$
	for any geodesic loop based at $x$ with length $\le \kappa l$.
\end{definition}
\subsection{Glossary of Main Results}
\subsubsection{Topological Remarks}
\begin{theorem}\label{thm:key-char}
	Suppose $\E^r$ is {\sf top-AC} with link $\N$ such that $\widetilde{\N}$ is $k$-connected; suppose in addition $n-1  \le (r-1)\wedge k$ holds. Then
	\begin{enumerate}
		\item  If $\M^n$ is closed, then $\widetilde{\M}$ is the $n$-sphere; the link is -- up to homotopy equivalence -- a $\sph^{r-1}$-bundles over $\sph^{n}$.
		\item If $\M^n$ is open, then it must be a contractible open manifold; in this case, $\N$ is a $(n+r-1)$-sphere and $\E \underset{\text{\sf diffeo}}{\cong} \R^{n+r}$.
	\end{enumerate}
\end{theorem}
\begin{corollary}\label{cor:top-ale-open}
	Suppose $\E^r$ is a {\sf top-ALE} vector bundle over an open manifold $\M^n$, then $\M$ is contractible consequently $\E$ is trivial hence a {\sf top-AE}. 
\end{corollary}
\begin{theorem}\label{thm:top-ale-general}
	Suppose $\E^r$ is a {\sf top-ALE} vector bundle over a closed $\M^n$, then
	\begin{enumerate}
		\item $\widetilde{\E}_{\sf slit}$ has a nonvanishing Euler class; consequently so does ${\E}_{\sf slit}$ given $\E$ is an orientable bundle;
		\item if $r\ge n$, then $\widetilde{\M}=\sph^n$ and $n=r$;
		\begin{enumerate}
			\item if in addition $\pi_1(\M)=\pi_1(\N)$, then $n=r = 2,4,8$; 
			\item  if in addition, the link of $\widebar{\E} = \widebar{\sf p}^*\E$ is an $H$-space for some normal covering $\widebar{\sf p}:\widebar{\M}\to \M$, then $n=r=4$;
		\end{enumerate}
	\end{enumerate}
\end{theorem}
\begin{corollary}\label{cor:ale-char}
	Suppose $\TM^n$ is {\sf top-ALE} then,
	\begin{enumerate}
		\item If $\M$ is closed, $n=2$ and $\widetilde{\M} =\sph^2$ and $\widetilde{\E}$ has link $\mathbb{RP}^2$;
		\item If $\M$ is open, then $\M$ is an open contractible manifold thus, $\mathcal{T}M \underset{\text{\sf homeo}}{\cong}  \R^{2n}$ ( $\underset{\text{\sf diffeo}}{\cong} $ for $n\neq 2$).
	\end{enumerate}
\end{corollary}
\begin{theorem}\label{thm:spin-4}
	Suppose $\E^2$ is a proper {\sf top-ALE} vector bundle with vanishing second Steifel-Whitney class over an open orientable spin 4-manifold $\M$, then either $\M$ is simply connected or its fundamental group is of the types presented in Proposition~\ref{prop:spin}. In particular, there exists no spin $\E^2$ over a spin $\M^4$ with even fundamental group. 
\end{theorem}
\subsubsection{Geometric Remarks}
\begin{theorem}[Nonnegative Curvature]\label{thm:nonneg-curv}
	\hfill
	\begin{enumerate}
		\item Suppose $\left(\M^n,g\right)$ is open complete and {\sf top-ALE} with $sec\ge 0$ and a nontrivial soul of codimension $r$ with $\lfloor \nicefrac{n}{2} \rfloor + 1 \le r\le n-2 $, then $n$ is even and $=2r$ and the soul is homeomorphic to a quotient of $\sph^{r}$. 
		\item Suppose $\E^r$ is {\sf top-AC} over a closed $\M^n$ and $\widetilde{\N}$ is $k$-connected with $n-1  \le (r-1)\wedge k$. If $\E$ or $\M$ admits a complete nonnegatively Ricci curved metrics with asymptotic growth $v_0 > \nicefrac{1}{m}$, then $\M = \nicefrac{\sph^n}{\Upgamma}$ with $|\Upgamma|< m$;
		\item 	Let $\E^{r\ge 3}$ be a {\sf top-AC} bundle over an open $\M^3$, $\widetilde{\N}$ is $k$-connected for $k\ge 2$ and suppose that $\M$ admits a complete metric with nonnegative scalar curvature.  Then $\M \underset{\text{\sf diffeo}}{\cong} \R^3$. 
	\end{enumerate}
\end{theorem}
\begin{theorem}[Asymptotic Flatness]\label{thm:AF-to-AC}
	Suppose $\left(\E^r, g\right)$ with $r\ge n$
	\begin{enumerate}
		\item 	is weakly asymptotically flat and with unique tangent metric cone at infinity (in particular, if it satisfies strong asymptotic flatness) and $\E_{\sf slit}$ has vanishing $\pi_i$ for $2\le i\le n-1$;
		\item[] or, 
		\item is strongly asymptotically flat and with strong holonomy control at scale $\kappa \in (0,\nicefrac{1}{2})$ and with $b_1(\E_{\sf slit})=0$ (vanishing first Betti number);
	\end{enumerate}
	then,  $M$ is either a differential spherical space form or an open contractible manifold. In case (2), $\E$ is also {\sf top-ALE}. 
\end{theorem}
\par The following results are about special metrics on bundles constructed from a metric $g$ on $\M$; see~\textsection\thinspace\ref{sec:esp-mets} for the precise definitions of the following metrics.
\begin{itemize}
	\item \emph{``connection metrics'':} a metric $g_{{\sf cm}}$ on a bundle $\E$ constructed from a given fiberwise metric $g_\E$, a connection $\nabla^\E$ compatible with $g_{_\E}$ and a rotationally symmetric metric $g_{_0}$ on $\R^r$) such that the bundle map becomes a Riemannian submersion with totally geodesic fibers $ \underset{\text{\sf isom}}{\cong} \left( \R^r,g_{_0} \right)$;
	\item \emph{``natural metrics'':} metrics $g_{{\sf nat}}$ on the tangent bundle that are decoupled and make the tangent bundle a Riemannian submersion;
	\item \emph{``admissible metrics'':} metrics $g_{{\sf ad}}$ on tangent bundle that are decoupled metrics and are defined using lifts of a given affine connection $\nabla$ and with a vertical part of special form~\eqref{eq:ad-mets}.
\end{itemize}
\begin{theorem}[Special Metrics]\label{thm:special-mets}
	\hfill
	\begin{enumerate}
		\item if $\left(\E,  g_{\sf cm}\right)$  is {\sf Riem-ALE}, then $\left(\M,g\right)$ is flat and $\nabla^\E$ is a flat affine bundle connection.
	\item if $\left(\TM,  g_{\sf nat}\right)$ is {\sf Riem-ALE}, then $\left(\M,g\right)$ is flat;
		\item  if $\left(\TM,  g_{\sf ad}\right)$ is  {\sf Riem-ALE}, then $\nabla$ is a flat affine connection. 
	\end{enumerate}
	In all above cases, $\M$ is open contractible and if furthermore $\M$ is $\nabla$-geodesically complete (in the first two cases $\nabla$ is the Levi-Civita connection), then $\M \underset{\text{\sf diffeo}}{\cong} \R^n$.  
\end{theorem}
\begin{remark}\label{rem:open-cont}
	There are examples of composite Whitehead manifolds $\widetilde{\M}$ that infinitely cover an open manifold $\M$. As such $\widetilde{\E}$ is {\sf top-AE} yet $\E$ is not {\sf top-ALE}. In particular, $\pi_1(\M)$ has no torsion or subgroups with finite index~\cite{TW}.
\end{remark}
\addtocontents{toc}{\protect\setcounter{tocdepth}{-1}}
\section*{\small \bf  Acknowledgements}
\addtocontents{toc}{\protect\setcounter{tocdepth}{1}}
\vspace{-165pt}
\begin{minipage}[t][7cm][b]{0.87\textwidth}
	\footnotesize
	\renewcommand{\labelitemi}{\raisebox{2pt}{\scalebox{1.2}{$\centerdot$}}}
	\begin{itemize}
		\item FA was partially supported by the Ministry of Science funding to visit IUT.
		\smallskip
		\item SL acknowledges partial support by the \emph{Kazemi Ashtiani early career award} 	and partial support from the resident researcher grant from \emph{IPM, Grant No. 1402530316}. 
		\smallskip
		\item ZF is in parts supported by the non-resident researcher grant from \emph{IPM, Grant No. 1402530045}.
	\end{itemize}
\end{minipage}
\normalsize
\vspace{10pt}
\section{Topological Rigidity}
Vector bundles over $\M$ suffer less topological (or geometric) anomalies than $\M$ does. Indeed, all bundles with rank $\ge 2$ over closed manifold admit {\sf top-AC} end structure; Even knowing the {\sf Riem-AC} end info of $\E$, one cannot in general pinpoint the topological and differential structure of $\M$. To demonstrate this, we start with a look at the some more or less well-known examples about the tangent bundles, as prototypes of vector bundles.
\subsection{Wildness of the Preimage under the Functor $\mathcal{T}$}
\hfill
\begin{example}\phantom (different homeomorphism types with diffeomorphic Euclidean tangent bundles) 
	Suppose $\mathcal{W}^n$ ($n\ge 3$) is any open contractible manifold, then $\mathcal{T}\mathcal{W}$ is a trivial bundles hence must be homeomorphic to $\R^{2n}$~\cite{McM}, thus also diffeomorphic to $\R^{2n}$.  By~\cite{McM2,CK,Glaser}, we know there exist uncountably many topologically different open contractible manifolds $\mathcal{W}^n$ for $n\ge 3$. 
\end{example}
\begin{example}\phantom (homeomorphic and not diffeomorphic with the same Euclidean tangent bundle)
	An exotic $\R^4$ whose tangent bundle is again trivial and diffeomorphic to $\R^8$.
\end{example}
\begin{example} \phantom (Lens spaces) Examples of non-homeomorphic closed manifolds with diffeomorphic tangent bundles include any two 3D Lens spaces $L(p,q_1)$ and $L(p,q_2)$ with 
	\[
	\pm q_1q_2 \equiv n^2\tpmod p, \quad q_1 \not\equiv \pm q_2^{\pm 1} \tpmod p.
	\]
	These have the same homotopy type~\cite{Whitehead}, hence their product with $\R^3$ (recall all three manifolds are parallelizable) are diffeomorphic~\cite{Mil}, yet they are not homeomorphic~\cite{Brody}.
\end{example}
\subsection{Homotopy Characterizations}
\begin{lemma}\label{lem:key}
		\hfill
	\begin{enumerate}
		\item Suppose $\E^{r\ge 2}$ (over the manifold $\M^n$) is {\sf top-AC} with end $\R^+\times \N$, then $\pi_m(\M) \le \pi_m\left(\N\right)$ for all $1 \le m < r$;
		\item $\E^{r\ge2}_{\sf slit}$ deformation retracts onto $\N$; hence the link is unique up to homotopy equivalence;
		\item Every vector bundle $\E^{r\ge 2}$ over a compact $\M^n$ is regular {\sf top-AC} with link that is unique up to homotopy equivalence and a $\sph^{r-1}$-bundle over $\M^n$.
		\end{enumerate}
\end{lemma}
\begin{proof}[\footnotesize \textbf{Proof}]
	\hfill
	\begin{enumerate}
		\item[(1)]
	Let $f:\sph^m \to \E$ be a continuous map based at the point $p=(x,q)$ in $\E\setminus \K$ where $q\in \N$. It is sufficient to base-homotope $f$ to a map with image included in $\{ q \}\times \N$. The proof consists of four steps:
	\begin{enumerate}
		\item[(i)] (Smooth approximation) We note that to compute homotopy groups in manifolds, by approximation theorem, we can work with smooth homotopies of smooth pointed maps from $\sph^m$ into $\E$; e.g. see~\cite[Corollary 17.8.1]{BT}. Therefore, we assume $f$ is smooth and  $f(a) = p$ for a base point $a\in \sph^m$. 
		\item[(ii)] (Local transversalization) Since 
		\[
		\dim(\M)+ m = m + n < n + r = \dim(\E),
		\]
		and since there exists an open neighborhood $\mathcal{U}$ of $a\in \sph^m$ and an open neighborhood $\mathcal{W}$ of $\M$ so that $f(U) \subset \E \setminus \mathcal{W}$, by the transversality theorem (e.g. \cite[Theorem 2.1]{Hirsch}), we can smoothly base-homotope $f$ to a map $\widebar{f}$ with $\widebar{f} \transv_{\sph^m} \M$ (i.e. that is tranversal to $\M$ along $\sph^m$). Considering the dimension deficiency this means the image of $\widebar{f}$ is disjoint from $\M$ (this is just the usual puling apart geometrically nontransversal objects in differential topology). For more details on this transversalization away from an open set, see~\cite[Corollary 4.1.2(b)]{GG}.
		\item[(iii)] (Cut-off dilation) Consider an open neighborhood $\mathcal{V}$ of $a\in \sph^m$ with $f(\mathcal{V}) \subset \E \setminus \K$. Let $0 \le \varphi \le 1$ be a smooth cut-off function on $\sph^m$ with 
		$
		\varphi(a) = 1
		$
		and
		$ \supp(\varphi) \subset \mathcal{V}$.
		\par Let $\Psi_\chi^t$ denote the flow of the Euler vector field $\chi$ of $\E$, then
		\[
		\lambda \mapsto f_\lambda := \Psi_\chi^{\lambda (1 - \varphi(\cdot) )} \left(  f(\cdot) \right) : \sph^m \to \E,
		\]
		is a smooth homotopy of $f$ that keeps the image of $V$ fixed; in particular, is a based homotopy. Since the image of $f$ does not intersect $\M$, from compactness of $\K$, one deduces $\Ima\left( f_\lambda \right) \subset \E\setminus \K$ for sufficiently large values of $\lambda$.
		\item[(iv)] (Radial retraction) Now since $\N$ is a deformation retraction of $\E \setminus \K$, by shrinking the radial cone fibers, we get a map with image included in $\N$.  
	\end{enumerate}
	Steps I-III provide us with the desired based homotopy. This shows $\pi_{m}\left( \E \right) \le \pi_{m}\left( \N \right) $ for all $m<r$. 
	\item[(2)] The proof follows from the dilation and radial retraction similar to item (1).
	\item [(3)]	It is standard that one can choose a fiberwise metric on $\E$. Then, $\E_{\sf slit}$ is diffeomorphic to $\R^+\times \N$ where $\N$ is the unit sphere bundle with respect to this fiberwise metric. $\E$ is even regular since taking $\K$ to be a disk bundle with respect to the chosen fiber-wise metric, $\E \setminus \K$ is again a topological punctured cone with link $\N$ and this time the homeomorphism extends to $\partial \K$.
\end{enumerate}
\end{proof}
		\begin{proposition}[periodicity]
		Suppose $\E^{r\ge 2}$ is a {\sf top-AC} whose link $\N$ has the same higher homotopy groups (from $\pi_5$ onward) as $\sph^{\ell}$. Then,
		\[
		\pi_{\ell+4}\left( \M^n \right) = \pi_{\ell+4}\left( \sph^{r-1} \right), \quad \text{whenever} \quad \ell \ge 7;
		\]
		in particular, for all {\sf top-ALE} bundles $\E^{r\ge 2}$,
		\[
		\pi_{n+r+3}\left( \M^n \right) = \pi_{n+r+3}\left( \sph^{r-1} \right), \quad \text{whenever} \quad n+r \ge 8.
		\]
	\end{proposition}
		\begin{proof}[\footnotesize \textbf{Proof}]
			Suppose $\pi_m(\N)=\pi_{m+1}(\N) = 0$. By the long exact sequence of homotopies, for the fibration ${\R^r}^*\hookrightarrow \E_{\sf slit} \to \M$, we get the long exact sequence
			\[
			\cdots \to \pi_{m+1}\left(\E_{\sf slit} \right) \to	\pi_m\left( {\R^r}^* \right) \to \pi_m(\M) \to \pi_m\left( \E_{\sf slit} \right) \to \cdots.
			\]
			By Lemma~\ref{lem:key}, $\pi_{k}\left(  \E_{\sf slit} \right) = \pi_k\left( \N \right)$ holds for all $k$; in particular $\pi_{m+1}\left(\E_{\sf slit} \right) = \pi_{m}\left(\E_{\sf slit} \right) = 0$. This implies $\pi_m(\M) = \pi_m({\R^r}^*)= \pi_m(\sph^{r-1})$. 
			
Now using the Frudehnthal's suspension theorem and low dimensional sphere homotopy groups (or Frudehnthal's suspension theorem combined with the list of stable homotopy groups; e.g. see~\cite{IWX}), we know $\pi_{l+4}(\sph^\ell) = \pi_{l+5}(\sph^\ell) = 0$. 
	\end{proof}
\subsection*{\small Proof of Theorem~\ref{thm:key-char}}
\par Suppose $\M^n$ is closed. By Lemma~\ref{lem:key}, we know $\pi_1(\M)$ is finite and $\widetilde{M}$ is $(n-1)$-connected. By the Hurewicz's theorem, $\pi_n(\widetilde{\M}) \equiv H_n(\widetilde{\M}) = \Z$. Now consider the continuous map $f:\sph^n \to \widetilde{\M}$ inducing the $n$-th homotopy class. This map clearly induces isomorphism on all homology groups hence by Hurewicz's theorem (needs simply connectedness), it is a weak homotopy equivalence hence by Whitehead's theorem, it is also a homotopy equivalence. Therefore, $\widetilde{\M}$ is a topological $n$-sphere. In particular, $\widetilde{\E}^{r\ge 2}$ is a {\sf top-AC} vector bundle over $\sph^n$ and by Lemma~\ref{lem:key}, the link, up to homotopy, is a $\sph^{r-1}$-bundle over $\sph^n$. 
\par If $\widetilde{\M}$ is open, then by $(n-1)$-contentedness and openness, one deduces that $\widetilde{\M}$ has trivial homology hence again by Hurewicz's theorem, $\widetilde{\M}$ also has trivial homotopy groups $\pi_i\left(\widetilde{\M}\right)$ for $i=1,\cdots,n$. So again, by Whitehead's theorem, the inclusion map of a point in $\widetilde{\M}$ is a homotopy equivalence i.e. $\widetilde{\M}$ is contractible. By~\cite[Lemma 2.2]{Wright}, an open contractible manifold cannot non-trivially and finitely cover another manifold. Therefore, $\M$ is itself contractible. 
\par As a consequence $\E^{r\ge 2}$ is diffeomorphic to $\M \times \R^r$ and by \cite{McM}, this is homeomorphic (and also diffeomorphic if $n+r\neq 4$) to $\R^{n+r}$. So, the link is a homotopy sphere thus a topological sphere (by the argument from the previous paragraph) and it follows form Lemma~\ref{lem:key} that the link of $\E^{r\ge 2}$ is homeomorphic to $\sph^{n+r-1}$.\qed
\par For even dimensional closed $\M$, we get a stronger conclusion.
\begin{theorem}
	Suppose $\E^{r\ge 2}$ is a vector bundle over closed $\M^{2n}$ with link $\N$ such that $\widetilde{\N}$ is $k$-connected with $n  \le (r-1)\wedge k$. Then
	$\widetilde{\M}$ is the $2n$-sphere.
\end{theorem}
\begin{proof}[\footnotesize \textbf{Proof}]
	\par A corollary of Whitehead's theorem is the minimal cell structure for manifolds with given finitely generated homology groups; e.g. see \cite[Proposition 4.C.1]{Hatcher}. As a famous consequence, an $(n-1)$-connected closed simply connected manifold $\widetilde{M}^{2n}$ has the homotopy type of a $2n$-dimensional disk attached to the wedge of $b_n$ number of $2$-spheres via an attaching map $f$ where $b_n$ is the $n$-th Betti number $b_n = \mathrm{rank}(H_n)$; see~\cite{Mil2} for an overview.
	\par So by the hypothesis, we deduce $\widetilde{\M}^{2n}$ is an $n$-connected hence has the homotopy type of an $2n$-disk with its boundary attached to a point; hence $\widetilde{\M}$ has the homotopy type of the $2n$-sphere therefore is homeomorphic to $\sph^{2n}$.
\end{proof}
\subsection{Coverings}
\par Suppose $\E^{r\ge 2}$ is {\sf top-AC} with link $\N$ and $\widebar{\sf p}:\widebar{\M}\to \M$ is a covering. Then one can equivariantly lift $\E$ to a vector bundle $\widebar{\E}^{r\ge 2}$ over $\widebar{M}$ i.e. $\widebar{\E}= \widebar{\sf p}^*\E$.
\begin{proposition}\label{prop:ac-cover}
	Let $\widebar{\sf p}:\widebar{\M}\to \M$ be a normal covering
	\hfill
	\begin{enumerate}
		 \item if $\M$ is a simple {\sf top-AC}. Then $\widebar{\M}$ is simple $\sigma$-{\sf top-AC}. If the covering is finite, then each connected component $\widebar{\mathcal{C}}_{\alpha j}$ of $\widebar{\sf p}^{-1}(\M\setminus \K)$ is a simple {\sf top-AC} end, $\nicefrac{\R^+\times \N_j}{\widebar{\Upgamma}_{\alpha j}}$ that covers the {\sf top-AC} end $\mathcal{M}_j = \nicefrac{\R^+\times \N_j}{\Upgamma_{j}}$ of $\M$ and we have $\widebar{\Upgamma}_{\alpha j}  \mathrel{\unlhd}  \Upgamma_j$. 
		 \item if $\E$ is a vector bundle over a compact $\M$ and $\widebar{\sf p}$ is a finite covering, then $\widebar{\N}$ as a $\sph^{r-1}$ fiberation over $\widebar{\M}$ is -- up to homotopy equivalence -- the lift (or pull-back bundle) of $\N$ via the covering $\widebar{\sf p}$; in particular, the deck group of $\widebar{\sf p}$ satisfies $\widebar{\G} = \nicefrac{\pi_1^\infty(\E)}{\pi_1^\infty(\widebar{\E})}$. 
		\item In general, if $\widebar{\sf p}$ is a finite covering then, $\E$ is a simple {\sf top-AC}  if and only if $\widebar{\E}$ is also a simple {\sf top-AC}; furthermore, $\widebar{\G} = \nicefrac{\pi_1^\infty(\E)}{\pi_1^\infty(\widebar{\E})}$ holds.
		\item Consequently, assuming $\pi_1(\M)$ is finite, $\E$ is a simple {\sf top-AC}  if and only if $\widetilde{\E}$ is also a simple {\sf top-AC} and $\pi_1(\M) = \nicefrac{\pi_1^\infty(\E)}{\pi_1^\infty(\widetilde{\E})}$ holds.
	\end{enumerate}
	Note that, working in the smooth category, the same statements hold with  {\sf diff-AC} replacing {\sf top-AC}.
\end{proposition}
\begin{proof}[\footnotesize \textbf{Proof}]
	\hfill
	\begin{enumerate}
		\item[(1)] Suppose $\K\Subset \M$ and $\mathcal{\M}_j$ is a {\sf top-AC} end. Set $\widebar{\K}:= \widebar{\sf p}^{-1}(\K)$. By the path lifting property, the uniqueness of lifts and the definition of normal coverings, one deduces that every path connected component $\widebar{\M}_{\alpha j}$ of $\widebar{\sf p}^{-1}\left(\M_j\right)$ is a normal covering of $\M_j$ with a deck group $\widebar{\G}_j^\alpha \mathrel{\unlhd} \widebar{\G}$. 
		Since the fundamental group of any manifold is countable, we deduce $\widebar{\K}$ is $\sigma$-compact.
	Consequently, 
			\[
		\widebar{G}_{j}^\alpha = \nicefrac{\pi_1\left(\mathcal{\M}_j\right)}{\widebar{\sf p}_* \pi_1\left(\widebar{\M}_j^\alpha\right)}, \quad
		\text{and} \quad  \G=\nicefrac{\pi_1(\M)}{\widebar{\sf p}_*\; \pi_1\left(\widebar{\M}\right)}.
		\]
		Any normal covering of the end  $\nicefrac{\R^+\times \N_j}{\Upgamma_{j}}$ is itself covered by the universal cover that is $\R^+\times \N_j$ thus  $\widebar{\M}_{\alpha j}=\nicefrac{\R^+\times \N_j}{\widebar{\Upgamma}_{\alpha j}}$ and $\widebar{\Upgamma}_{\alpha j}  \mathrel{\unlhd}  \Upgamma_j$ follows by standard normal covering properties; e.g. see~\cite[Proposition 1.40]{Hatcher} for further details. 
	\par When the covering is finite, $\widebar{\K}$ is also compact hence $\widebar{\M}$ becomes a simple {\sf top-AC}. 
	\item[(2)] Choosing a fiberwise metric on $\E$ and its lift on $\widebar{\E}$, one observes that the link $\widebar{\N}$ that is the fiberwise sphere bundle w.r.t to the chosen metric is also the lift of $\N$ under the covering map $\widebar{\sf p}$. This means $\N$ can be identified with the orbit space $\nicefrac{\widebar{\N}}{\widebar{\G}}$ and we have the normal covering $\widebar{\N} \to \nicefrac{\widebar{\N}}{\widebar{\G}}$. This readily implies the claim $\widebar{\G} = \nicefrac{\pi_1(\N)}{\pi_1(\widebar{\N})} = \nicefrac{\pi_1^\infty(\E)}{\pi_1^\infty(\widebar{\E})}$.
	\item[(3)] This is a straightforward consequence of item (1) and the fact that $\E^r$ and $\widebar{\E}^r$ both have only one ends along with uniqueness of link up to homotopy. 
	\item[(4)] a direct consequence of 2 and 3. 
	\end{enumerate}
\end{proof}
\subsection{\textbf{\sf Top-ALE} Vector Bundles}
\subsection*{\small Proof of Corollary~\ref{cor:top-ale-open}}
This is a direct consequence of Theorem~\ref{thm:key-char}. 
\subsection*{\small Proof of Theorem~\ref{thm:top-ale-general}}
	\begin{enumerate}
\item Suppose $\E$ is {\sf top-ALE}. By Proposition~\ref{prop:ac-cover}, $\widetilde{\E}$ is also {\sf top-ALE} and is an orientable bundle. Now suppose the Euler class of $\widetilde{\E}_{\sf slit}$ vanishes. There exists a closed global angular form and as a result, the Leray-Hirsch theorem applies~\cite[Chapter 11]{BT} hence the de Rham cohomology groups satisfy
 \[
 \mathsf{H}^*\left( \widetilde{\E}_{\sf slit}\right) = \mathsf{H}^*(\M) \otimes \mathsf{H}^*\left( \sph^{r-1}\right) =  \mathsf{H}^*\left(\M \times \sph^{r-1}\right),
 \]
as bi-graded algebras that is a contradiction e.g. by computing the $*=r-1$-th cohomology.
\item 
\begin{enumerate}
\item By Lemma~\ref{lem:key}, we deduce the link $\nicefrac{\sph^{n+r-1}}{\widetilde{\Upgamma}}$ of $\widetilde{\E}$ is a $\sph^{r-1}$-bundle over $\sph^n$. The exact Gysin sequence of cohomologies gives
\[
\cdots \to \mathsf{H}^{m}(\nicefrac{\sph^{n+r-1}}{\widetilde{\Upgamma}}) \to \mathsf{H}^{m-r+1}(\sph^n) \to \mathsf{H}^{m+1}(\sph^n) \to \mathsf{H}^{m+1}(\nicefrac{\sph^{n+r-1}}{\widetilde{\Upgamma}})\to \cdots;
\]
since $\left(\nicefrac{}{\widetilde{\Upgamma}}\right)^*$ is injective on cohomology groups, for $m=n-1$, we deduce $\mathsf{H}^{n-r}(\sph^n) = \mathsf{H}^{n}(\sph^n) = \Z$ thus $n=r$.  
\par Now if in addition, $\pi_1(\M)=\pi_1(\N)$, Proposition~\ref{prop:ac-cover} implies the link of $\widetilde{\E}$ is the sphere $\sph^{2n-1}$ and is a $\sph^{n-1}$-bundle over $\sph^n$. It can be argued that the Hopf invariant of the resulting bundle map ${\frak p}:\sph^{2n-1}\to \sph^{n}$ is $\pm 1$~\cite{Price}. For the sake of clarity, the sketch proof idea is the following.  The mapping cylinder $M_{\frak p}$ of ${\frak p}$ is a disk-bundle over $\sph^n$. Take $\tau$ to be the Thom class of this disk-bundle. By the relative Leray-Hirsch theorem, ${\frak p}^{*}\tau$ is the generator of  $\mathsf{H}^n({C_{\frak p}})$ where $C_{\frak {\frak p}}= \nicefrac{M_{\frak p}}{\sph^{2n-1}}$ is the mapping cone of $p$. Then by direct application of the Tom isomorphism, we deduce that the Hopf invariant is $\pm 1$; see~\cite{Price} for further details.  
\par So, after perhaps composing with an orientation reversing bundle isomorphism, the Hopf invariant of ${\frak p}$ is $1$; now by the classification in~\cite{Adams}, we deduce the claim on the possible dimensions. 
\item
\par The second claim directly follows item (2) and \cite{DS}. 
\end{enumerate}
\end{enumerate}

\begin{proposition}\label{prop:24}
	\hfill
	\begin{enumerate}
		\item {\sf Top-ALE} bundles on $\sph^2$ are in one-to-one correspondence with Lens spaces $\mathrm{L}_{m,1}$;
		\item  there exist 1 standard {\sf diff-AE} and 15 exotic {\sf diff-AE} bundles $\E^4$ over $\sph^4$.
	\end{enumerate}
\end{proposition}
\begin{proof}[\footnotesize \textbf{Proof}]
	\hfill
	\begin{enumerate}
		\item 	\par All circle bundles over the sphere are Seifert fiber manifolds. Consequently all nontrivial bundles are Lens spaces $\mathrm{L}_{m,1}$ where $m$ is the Euler number; furthermore,  they all come from vector bundles (meaning their structure group is $\mathrm{O}(n)$) showing that all non-trivial vector bundles over $\sph^2$ are {\sf top-ALE}. 
		\item By \cite{EK}, we know there exist 16 topological $\sph^7$ that have distinct differentiable structure and are $\sph^3$-bundles over $\sph^4$. These bundles are all extendable to rank $4$ vector bundles over $\sph^4$ since their structure group is $\mathrm{O}(4)$ and the resulting vector bundles must posses different differentiable structures (since their slit bundles obviously do).
	\end{enumerate}
\end{proof}
\begin{remark}
	A complete classification of $\sph^3$-bundles over $\sph^4$ up to homotopy, homeomorphism and diffeomorphism is carried out in~\cite{CE}.
\end{remark}
\begin{corollary}
	\hfill
	\begin{enumerate}
		\item There exists no {\sf top-ALE} bundle $\E^r$ over closed $\M^n$ with a non-vanishing global section;
		\item $\mathcal{T}\sph^{2n-1}$ is not {\sf top-ALE}.
	\end{enumerate}
\end{corollary}
\subsubsection{Classification of {\sf ALE} Tangent Bundles}
Aside from the Eguchi-Hanson spaces $\mathcal{T}\sph^2$ and the tangent bundle of open contractible manifolds, no other manifold admits {\sf ALE} tangent bundle. This is argued in this section. 
\begin{lemma}\label{lem:key-ALE}
	\hfill
	\begin{enumerate}
		\item 	Any asymptotic cone link for $\mathcal{T}\sph^n$ is homotopy equivalent to the Steifel manifold  \newline $\nicefrac{\mathrm{SO}(n+1)}{\mathrm{SO}(n-1)}$;
		\item For $n> 2$,	$\nicefrac{\mathrm{SO}(n+1)}{\mathrm{SO}(n-1)}$ and $\sph^{2n-1} = \nicefrac{\mathrm{SO}(2n)}{\mathrm{SO}(2n-1)}$ do not have the same higher homotopy groups. For $n=2$, the latter is a double cover of the former. 
	\end{enumerate}
\end{lemma}
\begin{proof}[\footnotesize \textbf{Proof}]
	\hfill
	\begin{enumerate}
		\item[(1)] First note that	$\mathcal{T}\sph^n_{\sf slit}$ is homeomorphic to $\R^+\times \nicefrac{\mathrm{SO}(n+1)}{\mathrm{SO}(n-1)}$. 
		To see this, let $\mathbb{Q}^n \subset \mathbb{C}^{n+1}$ be the affine hyperquadric given by the equation $\sum\limits_{j=1}^{j=n+1} z_j^2 = 1$. Then there is an $\mathrm{SO}(n+1)$-equivariant diffeomorphism $\varphi: \mathcal{T}\sph^n \to \mathbb{Q}^n$ given by
		\[
		(x,u) \to \cosh(\|u\|)x + i \left(\nicefrac{\sinh(\|u\|)}{\|u\|}\right) u,
		\]  
		presented in~\cite{Szoke}. Note that $(x,u)$ lie on the anti-de-Sitter space in $\R^{2n+2}$ and $u$ is considered as a vector in $\R^{n+1}$ that is perpendicular to $x$ (see the details in \cite{Szoke}). The isometric action of $\mathrm{SO}(n+1)$ on $\sph^n$ can be extended to the unit sphere bundle (i.e. set $\|u\| = 1$) in a natural way. The stabilizer of an $(x,u)$ is clearly $\mathrm{SO}(n-1)$. Hence the unit tangent bundle is diffeomorphic to $ \nicefrac{\mathrm{SO}(n+1)}{\mathrm{SO}(n-1)}$ so the slit tangent bundle is diffeomorphic to $\R^+\times \nicefrac{\mathrm{SO}(n+1)}{\mathrm{SO}(n-1)}$.
		\par Now the conclusion follows from Lemma~\ref{lem:key}. 
		\item[(2)] $\nicefrac{\mathrm{SO}(n+1)}{\mathrm{SO}(n-1)}$ and $\sph^{2n-1} = \nicefrac{\mathrm{SO}(2n)}{\mathrm{SO}(2n-1)}$ are indeed the Steifel manifolds $\mathrm{O}(n+1,2)$ and $\mathrm{O}(n+1,1)$. $\mathrm{O}(n,k)$ fits into the long exact sequence of homotopies 
		\[
		\cdots \to	\pi_m\left( \mathrm{O}(n-1,k-1) \right) \to \pi_m\left( \mathrm{O}(n,k) \right) \to \pi_m\left( \sph^{n} \right) \to \pi_{m-1}\left( \mathrm{O}(n-1,k-1) \right) \to \cdots.
		\]
		One consequence is $\mathrm{O}(n,k)$ is $n-k-1$-connected thus the first claim follows. Since $\sph^3$ is a double cover of $\mathbb{RP}^3 = \mathrm{SO}(3)$ the second claim follows.   
	\end{enumerate}
\end{proof}
\subsection*{\small Proof of Corollary~\ref{cor:ale-char}}
\par Suppose $\M^n$ is closed. Since $\TM^n$ is {\sf ALE}, by Lemma~\ref{lem:key}, we know the fundamental group of $\M$ is finite and furthermore the link - up to homotopy - is $\nicefrac{\sph^{2n-1}}{\Upgamma}$.  By Theorem~\ref{thm:key-char}, $\widetilde{\M} = \sph^n$ and by Proposition~\ref{prop:ac-cover}, $\mathcal{T}\M$ is {\sf ALE} with link $\nicefrac{\sph^{2n-1}}{\widetilde{\Upgamma}}$ for some $\widetilde{\Upgamma} \mathrel{\unlhd}  \Upgamma $. By Lemma~\ref{lem:key-ALE}, we deduce $n=2$. 
	\par The claim for open $\M$ is a direct consequence of Theorem~\ref{thm:key-char}. \qed
\subsubsection{Low Dimensional Subrank Cases}
\par Characterizing sub-rank {\sf ALE} bundles (i.e. when $r<n$) is more delicate of a problem and these bundles are the habitat of interesting phenomena. For example, the first counter example to the ``\emph{positive mass theorem}'' for {\sf Riem-ALE} manifolds constructed in~\cite{LeBrun} is a line bundle over $\mathbb{CP}^1$. 
\par We close this section by the following two sub-rank cases one with closed and other with open base. 
\begin{theorem}
No {\sf top-ALE} bundle $\E^2$ exists over $\sph^n$ for $n>2$.
\end{theorem}
\begin{proof}[\footnotesize \textbf{Proof}]
By \cite{Steenrod}, we know the only $1$-sphere bundle over $\sph^n$ ($n>2$) is a product bundle therefore not homotopy equivalent to a quotient of $\sph^3$.
\end{proof}
We refer the reader to~\cite{Wolf} for classification of spherical space forms and the terminologies used. 
\begin{proposition}\label{prop:spin}
	\par Spherical space form groups in $\mathrm{SO}(6)$ that admit isomorphic lift to $\mathrm{Spin}(6)$ are cyclic groups of odd order (type I with $d=1$) or specific types of type I groups with odd order (corresponding to $d=3$) with presentation
	\[
	A^m=B^n= \mathbf{1}, \quad BAB^{-1} = A^r,
	\]
	\[
	n\equiv 3\!\!\!\! \pmod 6, \quad m\equiv 1 \!\!\!\! \pmod 2, \quad \left( n(r-1), m  \right) = 1, \quad r \not\equiv r^3 \equiv 1 \!\!\!\! \pmod m.
	\]
\end{proposition}
\begin{proof}
	\par First recall the vector representation $\uprho: \mathrm{Spin}(n) \to \mathrm{SO}(n)$ which gives $\mathrm{Spin}(n)$ as a double cover of $\mathrm{SO}(n)$. The tangent map of this representation is 
	$
	d\uprho: \mathfrak{spin}(n) \to \mathfrak{so}(n)
	$,
	given by $d\uprho(e_ie_j) = 2(e_i\wedge e_j)$ where $e_i\wedge e_j$ is the generator of rotations in the $e_i-e_j$ plane; e.g. see~\cite{Jost}. Hence, the lifts of the element $\mathrm{Rot}(2\pi\theta)\in \mathrm{SO}(n)$ are $\pm \exp(\theta \pi e_ie_j)\in \mathrm{Spin}(n)$.
	\par Following the presentation of~\cite{Wolf}, type I groups with $d=1$ are cyclic and generated by a matrix 
	\[
	A= \begin{pmatrix}  \mathrm{Rot}(1\pi/m) & 0 & 0 \\ 0 & \mathrm{Rot}(2\pi a/m)  & 0 \\ 0 & 0 & \mathrm{Rot}(2\pi b/m)   \end{pmatrix},
	\]
	with $(a,m)=1$ and $(b,m)=1$. The lifts are thus given by
	\[
	\widebar{A} = \mathrm{sign}^{\widebar{A}}_A \exp\left(  \nicefrac{\pi}{m}\; \widebar{e}_{12} +   \nicefrac{\pi a}{m} \; \widebar{e}_{34} +   \nicefrac{\pi b}{m} \; \widebar{e}_{56} \right), \quad \mathrm{sign}^{\widebar{A}}_A = \pm 1,
	\]
	\[
	\widebar{A}^m = \left(  \mathrm{sign}^{\widebar{A}}_A  \right)^m (-1)^{1+ a + b} \; \mathbf{1}.
	\]
	If $m$ is even, then $1+a+b$ must be odd hence $\widebar{A}^m = - \epsilon_A^m  \mathbf{1} = - \mathbf{1}$, so in this case, we do not get the same presentation relations. If $m$ is odd, then, upon setting $ \epsilon_A = (-1)^{1+ a + b}$, one gets a lift of $A$ that satisfy the same idempotency. 
	\par Type I groups with $d=3$ are given by generators $A$ and $B$ where $A$ is as in before with 
	$
	a=r, b=r^2
	$
	and 
	\[
	B= 
	\begin{pmatrix}  0 &I & 0 \\ 0 & 0 & I \\ \mathrm{Rot}(3l/n) & 0 & 0  \end{pmatrix}.
	\]
	with conditions
	\[
	n\equiv 0 \!\!\!\! \pmod 3, \quad \left( n(r-1), m  \right) = 1, \quad r \not\equiv r^3 \equiv 1 \!\!\!\! \pmod m.
	\]
		\par As block matrices, one has
	$
	B^3 = \mathrm{Rot}(3l/n) \mathbf{1}
	$;
	Set $\widebar{B} :=  \mathrm{sign}^{\widebar{B}}_B \widehat{B}$, where $\widehat{B}$ is the lift of $B$ with 
	\[
	\widehat{B}^3 =  \exp\left(\nicefrac{3\pi l}{n} \; e_1e_2 +  \nicefrac{3\pi l}{n} \; e_3e_4 + \nicefrac{3\pi l}{n} \; e_5e_6 \right).
	\]
 By $\widebar{B^3} = \widebar{B}^3$, one deduces
	\[
	\widebar{B}^3 =  \left(   \mathrm{sign}^{\widebar{B}}_B \right)^3 \exp\left(\nicefrac{3\pi l}{n} \; e_1e_2 +  \nicefrac{3\pi l}{n} \; e_3e_4 + \nicefrac{3\pi l}{n} \; e_5e_6 \right),
	\]
	and consequently
	\[
	\widebar{B}^n = \left( \widebar{B}^3  \right)^{\frac{n}{3}} =  \left(   \mathrm{sign}^{\widebar{B}}_B \right)^n (-1)^{3l} \mathbf{1}.
	\]
	\par We already know $3|n$. If $n=6k$, we get
	$
	\widebar{B}^n = (-1)^{3l} \mathbf{1} 
	$
	hence $l$ must be even that contradicts $(l,n/3)=1$. For $n=6k+3$, we get
	\[
	\widebar{B}^n = \left(   \mathrm{sign}^{\widebar{B}}_B \right)^3  (-1)^{3l} \mathbf{1} = \left(  \mathrm{sign}^{\widebar{B}}_B \;  (-1)^l \right)^3 \mathbf{1},
	\]
	hence upon setting 
	$
	\mathrm{sign}^{\widebar{B}}_B = (-1)^l,
	$
	one gets a lift with the same idempotency relation.
\par Let us check the relation $\widebar{B}\widebar{A}\widebar{B}^{-1} = \widebar{A}^r$. By the lifting, we know
	$
	\widebar{B}\widebar{A}\widebar{B}^{-1} = \delta \widebar{A}^r
	$.
	Raising sides to the $m$-th power (that is odd), we get
	$
	\delta^m = 1 \Longrightarrow \delta=1
	$.
	This means we get at least one spin structure and $m$ and $n$ are odd. 
\end{proof}
\subsection*{\small Proof of Theorem~\ref{thm:spin-4}}
Suppose $\E^2$ is a spin vector bundle over a spin $4$-manifold that is proper {\sf top-ALE} with link $\nicefrac{\sph^5}{\Upgamma}$. By naturality of Steifel-Whitney classes, it is straightforward to show $\E$ is spinable as a manifold. Consequently, so is the open submanifold $\E\setminus \K$ which is homeomorphic to $\nicefrac{{\R^6}^*}{\Upgamma}$. This means $\Upgamma \le \mathrm{SO}(6)$ is a spherical space form group admitting isometric lift to  $\mathrm{Spin}(6)$. By Proposition~\ref{prop:spin}, $\Upgamma$ must be one of the groups presented and in particular must be an odd group. \qed
\subsection{Topological Rigidity of $\K$}
We have seen that both $\E_{\sf slit} = \E \setminus \M$ and $\E \setminus \K$ deformation retract to $\N$. Below, we show under some mild regularity conditions on $\partial \K$, we have their boundaries have the same fundamental group i.e. $\pi_1(\K) = \pi_1(\M)$. Also we furthermore discuss some implied homological information on $\K$. 
\par Recall a topological submanifold $\mathcal{W}^{n-1}\subset \M$ is called locally flat whenever every $b\in \mathcal{B}$ admits a neighborhood $\mathcal{N}_b$ and a homeomorphism $h_b: \mathcal{N}_b \to \R^n$ such that $h_b\left( \mathcal{N}_b \cap \mathcal{W} \right) \subset \R^{n-1}$. In other words, $\mathcal{W}$ is locally flat in $\M$ if it admits a $\mathcal{C}^0$ atlas comprised of of "adapted $\mathcal{C}^0$ charts" from $\M$.
\begin{definition}\label{defn:reg-comp}
	We say $\K\subset \M$ is regular when $\partial \K$ is a two-sided locally flat $\mathcal{C}^0$-submanifold of $\M$ and is homeomorphic to $\N$. 
\end{definition}
\begin{theorem}\label{thm:comp-rig-pi}
	Suppose $\E^r$ is {\sf top-AC} with link $\N$ and the omitted regular compact set $\K$; then $\pi_1(\M) = \pi_1(\K)$. 
\end{theorem}
\begin{proof}[\footnotesize \textbf{Proof}]
\par It follows from the hypotheses that $\partial \K$ admits tubular neighborhood~\cite{MB}. We can consider open fattenning and thinning of $\K$ as follows. Let $\psi: \partial \K \times (-\epsilon, \epsilon) \to \mathcal{U}_\K$ be a tubular neighborhood of $\partial \K$. Define the open sets
	\[
	\mathcal{U}_{\alpha}:=  \begin{cases}  \K \cup \psi \left( \partial \K \times (-\epsilon, \upalpha)\right), \quad \text{for}\quad 0<\upalpha<\epsilon \quad \text{fattening} \\
		\K \smallsetminus \psi \left( \partial \K \times [\upbeta, \epsilon)\right), \quad \text{for}\quad -\epsilon <\upbeta<0 \quad \text{thinning}
	\end{cases}.
	\]	
\par For $\upalpha>0$ and $\upbeta <0$ sufficiently close to $0$, set $U_1 := \K_\upalpha$ and $U_2 := \E \setminus \widebar{\K_\upbeta}$ thus $U_1\cap U_2 = \partial \K \times (\alpha,\beta)$ that is homotopy equivalent to $\partial \K = \N$. 
\par Choosing a base point $q \in \partial \K = \N$, we note that every based loop $\upgamma$ in $U_1\cap U_2$ is based-homotope to a loop in $\N$. By Seifert-Van Kampen theorem, we deduce
\[
\pi_1(M) = \pi_1(\E) = \nicefrac{\pi_1\left( U_1  \right) \star \pi_1\left( U_2 \right)}{\iota_* \pi_1\left( U_1\cap U_2 \right)} = \nicefrac{\pi_1\left( \K_\alpha \right) \star \pi_1(\N)}{\iota_* \pi_1\left( \partial \K \right)} = \pi_1(\K).
\]
where $\iota$ signifies the inclusions of $U_1\cap U_2$ into $U_1$ and $U_2$. 
\end{proof}
\begin{corollary}
	Suppose $\E^r$ is a simple {\sf top-AC} with omitted regular compact set $\K$ and $\widetilde{p}: \widetilde{\M} \to \M$ is the universal covering, then $\widetilde{p}^{-1}(\K) = \widetilde{\K}$.
\end{corollary}
\begin{theorem}
	Suppose $\M$ is $(m,m+1)$-locally connected for some $m\ge 0$ and $\E^r$ is a {\sf top-AC} with link $\N$ with the regular omitted compact set $\K$. Then, ${\sf H}_m(\K) = 0$. 
\end{theorem}
\begin{proof}[\footnotesize \textbf{Proof}]
Applying the Mayer-Vietoris to the sets $U_i$, $i=1,2$ from the proof of Theorem~\ref{thm:comp-rig-pi}, and noting the homotopy equivalences $\E \sim \M$ and $U_1\cap U_2 \sim \partial K$, $U_1\sim \K$ and $U_2 \sim \N$, we get an exact sequence
\[
\cdots \to {\sf H}_{m+1}(\M) \to {\sf H}_{m}(\N) \to   {\sf H}_{m}( \K) \oplus  {\sf H}_{m}(\N) \to {\sf H}_m(\M) \to \cdots,
\]
so from hypothesis, the two middle spaces are isomorphic, yielding ${\sf H}_{m}( \K) = 0$. 
\end{proof}
\begin{corollary}
	Suppose $\widetilde{\N}$ is $k$-connected, then $\widetilde{\K}$ has trivial homology up to degree $(r-2)\wedge (k-1)$.
\end{corollary}
\begin{corollary}
	Suppose $\E^r$ is {\sf top-ALE} with omitted regular compact set $\K$, then $\widetilde{\K}$ has trivial homology up to degree $r-2$. 
\end{corollary}
\par Let us show that a ``regular''  {\sf top-AC} manifold $\M$ (see Definition~\ref{defn:regular}), $\K$ actually is a regular omitted compact set in the sense of Definition~\ref{defn:reg-comp}. 
\begin{proposition}
	Suppose $\M$ is {\sf top-AC} with omitted compact set $\K$. Then, $\M$ is a regular {\sf top-AC} if and only if $\K$ is a regular omitted compact set. 
\end{proposition}
\begin{proof}[\footnotesize \textbf{Proof}]
\par By hypothesis, there exists a homeomorphism $\varphi: \R^+ \times \N \to  \M\setminus \K$ that extends to $\partial \K$. Therefore, $\partial K$ is homeomorphic to $\N$ i.e. homeomorphic to the boundary of a smooth manifold with boundary. This implies that $\partial \K$ is collared within $\varphi\left(\R^+ \times \N\right)$ since it equals $\varphi(\{0\}\times \N)$. Also being homeomorphic to the boundary of a smooth manifold, $\partial K$ is locally flat. This means $\partial \K$ is collared within $\K$ as well~\cite{MB}. Putting these together, $\partial \K$ is two-sided and locally flat. 
\par Conversely, if $\partial \K$ is a regular omitted compact set, by~~\cite{MB}, it admits a tubular neighborhood, this together with $\partial \K \underset{\textsf{homeo}}{\cong} \N$ enables us to extend $\varphi$ to $\partial \K$ and indeed beyond by extending $\varphi$ by a collar of $\partial \K$. 
\end{proof}
\section{Geometric Rigidity}
In this section, we briefly touch upon the rigidity of {\sf top-AC} vector bundles under imposition of further restrictions on a given Riemannian metric $\mathdutchcal{g}$; in particular, these all apply to {\sf Riem-ALE} bundles. 
\subsection{Nonnegative Curvature}
\par Recall the asymptotic volume growth of a Ricci non-negatively curved metric $g$ is defined by
\[
v_0 := \lim_{r\to \infty} \nicefrac{\mathrm{vol}_g(B_r(p))}{\mathrm{vol}_{0}(B_r(0))} \in [0,1].
\]
\subsection*{\small Proof of Theorem~\ref{thm:nonneg-curv}}
\begin{enumerate}
	\item Let $\M$ be open and with nonnegative sectional curvature; then, by Cheeger-Gromoll's soul theorem, $\M$ is a normal bundle over a compact totally geodesic soul $\mathcal{S}$. Since the codimension of the soul is larger than half of the dimension, item (2) in Theorem~\ref{thm:top-ale-general} verifies the claim. 
	\item 
	By~\cite{Li, Anders}, we know the fundamental group of $\E$ (hence that of $\M$) is finite with less than $m$ elements. The claim then follows from the topological classification theorem. 
	\item Using the topological rigidity, we deduce $\M$ is open contractible 3 manifold. By the recent results in~\cite{Wang} regarding nonnegative scalar curvature, it is then diffeomorphic to $\R^3$. 
\end{enumerate}
\qed
\subsection{Asymptotic Flatness}
\subsection*{\small Proof of Theorem~\ref{thm:AF-to-AC}}
	\begin{enumerate}
		\item By \cite{PT}, we know $\E^r$ is {\sf AC}, then we apply Theorem~\ref{thm:key-char}. 
		\item By \cite{CL}, $\E$ is {\sf AC} and its end is a $\mathbb{T}^{b_1}$-fibration over an {\sf ALE} end. 
	\end{enumerate}
\qed
\subsection{Special Metrics on Bundles}\label{sec:esp-mets}
\begin{definition}
	A function $F:\R^m\to \R$ with $F(p)=0$ is admissible whenever there exists sub(super)-harmonic function $G$ with $G(p)=0$ such that $\lambda G \le F \le \Lambda G$ for $0< \lambda < \Lambda < \infty$. 
\end{definition}
It is easy to see that a Liouville-type theorem holds for admissible functions i.e. boundedness of $F$ implies $F \equiv 0$. 
\begin{lemma}\label{lem:sub-flat}
	Suppose $\mathdutchcal{g}$ is a Riemannian metric on $\E^r$ and $g$, is one on $\M$ with respect to which ${\frak p}: \E \to \M$ is a Riemannian submersion; let $\mathcal{V}_{(x,u)} = \mathcal{T}_{(x,u)}\E_{x}$ denote the vertical lifts. Suppose for every fixed $x$ and $X,Y\in \mathcal{T}_x\M$, the function
	\[
	F_{X,Y,x}(u) := \|\mathrm{Proj}_{\mathcal{V}_{(x,u)}}[X^h, Y^h]_{(x,u)} \|^2_{g_{_{\E_x}}}: \E_x \to \R,
	\]
	is admissible. If $\left(\E,\mathdutchcal{g}\right)$ is {\sf Riem-ALE} then $\left(\M, g\right)$ is flat. In particular, if $r\ge n$, then $\M$ is an open contractible manifold. 
\end{lemma}
\begin{proof}[\footnotesize \textbf{Proof}]
By O'Neil's fundamental equations for Riemannian submersions, we have
\[
 \|\mathrm{Proj}_{\mathcal{V}_{(x,u)}} [X^h, Y^h]_{(x,u)}\|^2 = \nicefrac{4}{3}\left( \mathrm{K}_{\E}(X,Y)_{(x,u)} - \mathrm{K}_\M(X,Y)_x \right).
\]
We observe that the RHS is bounded since $\E$ is assumed to be {\sf Riem-ALE}. The admissibility implies 
$
\mathrm{K}_{\E}(X,Y)_{(x,u)} = \mathrm{K}_\M(X,Y)_x
$.
Replacing $u$ by $\lambda u$ and letting $\lambda\to \infty$, we see that this quantity must be zero since $\mathrm{K}_{\E}(X,Y)_{(x,u)} $ vanishes at infinity by the {\sf Riem-ALE} assumption. 
\par When $r\ge n$, by topological rigidity and noting that that $\sph^n$ does not admit flat metric (e.g. by Cartan-Hadamard theorem), the contractability claim follows. 
\end{proof}
\subsubsection{Connection Metrics on $\E$}
Choosing a fiberwise metric $g_{{\E_p}}$ and a vector bundle connection $\nabla^\E$  compatible with $g_{{\E_p}}$, a metric $g$ on $\M$ and a rotationally symmetric metric $g_{_0} = dr^2 + \psi(r)^2g_{_{\sph^{r-1}}}$ on $\R^r$, there exists a unique metric $\mathdutchcal{g}_{{\sf cm}}$ (connection metric) with respect to which ${\frak p}: \E \to \M$ is a Riemannian submersion with horizontal distribution determined by $\nabla^\E$ and with totally geodesic fibers isometric to $\left(\R^r, g_{_0}\right)$; see~\cite{Besse}. 
\subsection*{\small Proof of Theorem~\ref{thm:special-mets}, Part I}
	By~\cite[Lemma 2.2]{SW}, we know for $u\neq 0$
	\[
	F_{X,Y,x}(u) = -\nicefrac{3}{4} \; \left(\nicefrac{\psi(\|u\|)}{\|u\|}\right)^2 \left| \mathrm{Rm}^{\nabla^\E} (X,Y)(u) \right|^2,
	\]
and vanishes for $u=0$. Since $\left| \mathrm{Rm}^{\nabla^\E} (X,Y)(u) \right|^2$ is sub-harmonic in $u$, we deduce $F_{X,Y,x}(u)$ is admissible  whenever $|\psi(r)| \le O(r)$ as $r\nearrow \infty$. By {\sf ALE} assumption the latter decay holds on $\psi$ thus, $\nabla^\E$ is flat affine since it has vanishing curvature and by Lemma~\ref{lem:sub-flat}, we deduce $\left(M,g\right)$ is flat. The second claim follows from topological rigidity.\qed
\subsubsection{Natural and Decoupled {\sf ALE} Metrics on $\TM$}
Given a linear connection $\nabla$ on $M$, any Riemannian metric $\mathdutchcal{g}$ on $\TM$ can be written in the form
\begin{align}
	\begin{cases}
		\mathdutchcal{g}_{(p,u)} \left( X^h, Y^h \right) =  g_{(p,u)}^h(X, Y)  \\
		\mathdutchcal{g}_{(p,u)} \left( X^h, Y^v \right) = a_{(p,u)}^{c}(X,Y)\\
		\mathdutchcal{g}_{(p,u)} \left( X^v, Y^v \right)  = g_{(p,u)}^v(X, Y)  
	\end{cases},
\end{align}
where $g^h_{(\cdot,u)}$ and $g^v_{(\cdot, u)}$ constitute families of Riemannian metrics on $\M$ parameterized by $u$ and $a^c_{(\cdot, u)}$ is a symmetric 2-tensor. Given connection coefficients $\mathcal{H}_i^j = u^b\Gamma_{ib}^j$ in a standard local coordinates $(x^i,u^j= \partial_j)$ of $\TM$, any metric $\widetilde{g}$ on $TM$ can be expressed in the form
\[
\mathdutchcal{g} = g^h_{ij} (\delta x^i)^\flat \otimes (\delta x^j)^\flat + a^c_{ij}  (\delta x^i)^\flat  \otimes du^j  + a^c_{ij}  du^i \otimes  (\delta x^j)^\flat + g^v_{ij} du^i \otimes dy^j,
\]
in which
$
\delta x^i := \nicefrac{\partial}{\partial x^i} - \mathcal{H}_i^j \nicefrac{\partial}{\partial u^j}
$.
\begin{definition}
	We say $\mathdutchcal{g}$ is decoupled, when $a^{c}(X,Y) \equiv 0$.
\end{definition}
The following is a mainstream notion. 
\begin{definition}
	We say $\mathdutchcal{g}$ is a natural metric whenever it is decoupled and furthermore, $g^h = g$ is independent of $u$ and $\nabla$ is the Levi-Civita connection of a the metric $g$. In other words, $\mathdutchcal{g}$ is natural whenever it is decoupled with respect to which $\boldsymbol{\pi}:\left(\TM, \mathdutchcal{g}\right)\to\left(\M, g\right)$ is a Riemannian submersion. We will denote natural metrics by $\mathdutchcal{g}_{\sf nat}$. 
\end{definition}
\subsection*{\small Proof of Theorem~\ref{thm:special-mets}, Part II}
This is a direct consequence of Lemma~\ref{lem:sub-flat} by noting that for a natural metric, the O'Neil's formula gives
\[
\mathrm{K}_{\E}(X^h,Y^h)_{(x,u)} = \mathrm{K}(X,Y)_{x} - \frac{3}{4} \left\|  (\mathrm{R}(X,Y) u)_x^v   \right\|^2,
\]
for orthonormal pair $X,Y$~\cite[Corollary 6.4]{GudKap}. \qed
\par Recall the Sasaki and Cheeger-Gromoll metrics~\cite{Sas,CG} are two natural metrics respectively determined by
$
\mathdutchcal{g}_{_{\sf S}}\left(X^v,Y^v\right) := g_p(X,Y)
$, 
and
$
\mathdutchcal{g}_{_{\sf CG}}(X^v,Y^v) := \nicefrac{ \left( g_p(X,Y) + g_p(X,u)g_p(Y,u)  \right)}{1+r^2}
$.
\begin{remark}
	By \cite{Aso}, flatness of $g_{_{\sf S}} $ would follow from much weaker assumption of boundedness of the sectional curvature of $\TM$.
\end{remark}
\begin{proposition}
	$\left(\TM, \mathdutchcal{g}_{_{\sf CG}} \right)$ is never {\sf Riem-ALE}. 
\end{proposition}
\begin{proof}
By \cite[Corollary 8.9]{GudKap}, 
	\[
	\widetilde{\mathrm{K}}(X^v, Y^v) = \nicefrac{(1-\alpha)}{\alpha^2} + \left(\nicefrac{\alpha+2}{\alpha}\right) \nicefrac{1}{\left(1 + g(X,u)^2 + g(Y,u)^2\right)}.
	\]
	Therefore, if we choose $u \perp \mathrm{span}\{X,Y\}$, taking the vectors $\lambda u$ and letting $\lambda \to \infty$, we get
	$
	\widetilde{\mathrm{K}}(X^v, Y^v)_{(p,\lambda u)} \to 1
	$,
which violates being {\sf Riem-ALE}. 
\end{proof}
Below we show flatness for another family of metrics on the tangent bundle not coming from Riemannian submersions. 
\begin{definition}[admissible metrics]
	Consider a coordinate chart $(\Omega, \psi)$ and the metric $g_{_0}:=\psi^*\delta$. 
	\begin{enumerate}
		\item 	A metric $\mathdutchcal{g}$ is said to be \emph{admissible} in the domain $\mathcal{T}\Omega$ of a coordinate chart whenever it is \emph{decoupled} and we have
		\begin{align}\label{eq:ad-mets}
		g^v_{(p,u)} = \varphi(p,u) {\sf g}_p + (1-\varphi) g_{_0}, \quad \varphi \in \mathcal{C}^\infty\left( \mathcal{T}\Omega, [0,1] \right),
		\end{align}
		for some Riemannian metric ${\sf g} \neq g_{_0}$ on $\Omega$. 
		\item Suppose $\TM$ admits standard coordinates $x^i, u^i$ at infinity (outside a compact set $\K$). A metric $\mathdutchcal{g}$ on $\TM$ is called ``admissible at infinity'' when it is admissible on $\TM \setminus \K$ with respect to coordinates $x^i,u^i$. Below $g_{\sf ad}$ denotes an admissible metric at infinity of $\TM$. 
	\end{enumerate}
\end{definition}
\subsection*{\small Proof of Theorem~\ref{thm:special-mets}, Part III}
Suppose $\left(\TM, \mathdutchcal{g} = \mathdutchcal{g}_{\sf ad} \right)$ is {\sf Riem-ALE} with standard coordinates $x^i, u^i$ at infinity. Using \eqref{eq:ad-mets}, 
and setting $ \partial_{\widebar i}:=	\partial_i^v$, we have $	\partial_i^h  = \partial_i - \mathcal{H}_i^k \partial_{\widebar k}$, we get the metric local entries
\[
\begin{cases}
	\mathdutchcal{g}_{ij} &= g^h_{ij}   + \mathcal{H}_j^l  a^{c}_{il} + \mathcal{H}_i^ka^{c}_{jk} + \mathcal{H}_i^k\mathcal{H}_j^l g^v_{kl}  \\[10pt]
	\mathdutchcal{g}_{i\widebar{j}} &=  \mathcal{H}_i^k \phi{\sf g}_{kj}  + (1-\varphi) \mathcal{H}_i^k
	\\[10pt]
	\mathdutchcal{g}_{\widebar{i}\widebar{j}} &= \varphi {\sf g}_{ij} + (1-\varphi)\delta_{ij}
\end{cases}.
\]
\par From the {\sf ALE} hypothesis, we get
$
\lim_{\|u\| \to \infty}\varphi(p,u) = 0
$. 
So from $\lim_{\|u\| \to \infty} \mathdutchcal{g}_{i\widebar{j}} =0 $ we deduce $\lim_{\|u\| \to \infty} \mathcal{H}_i^j = 0$. Now, for fixed $i,j$, and at fixed $p$, by letting
$
u^b = t\Upgamma^j_{ib}
$,
and using $\mathcal{H}_i^j = u^b\Gamma_{ib}^j$ , upon letting $t\nearrow \infty$, we deduce $\Upgamma^j_{ib}(p) = 0$. Hence $\nabla$ is a flat connection.
\qed
\begin{remark}
	\par Recall by classification of flat manifolds and by applying Cartan-Hadamard theorem, a complete contractible flat $n$-manifold is diffeomorphic to $\R^n$. Even more generally, using Cartan-Ambrose-Hicks theorem, we deduce a contractible flat affine manifolds (w.r.t. a given $\nabla$) that is $\nabla$-geodesically complete must be diffeomorphic to $\R^n$. 
\end{remark}
\renewcommand{\baselinestretch}{1.1}
\vspace{20mm}
\hfill \textsf{\scriptsize -- \today}

\end{document}